\documentclass[letterpaper, 10pt, conference]{ieeeconf}      

\IEEEoverridecommandlockouts                              
\newcommand{\ud}{\,\mathrm{d}}

\usepackage{graphics} 
\usepackage{epsfig} 
\usepackage{mathptmx} 
\usepackage{times} 
\usepackage{amsmath} 
\usepackage{amssymb}  
\usepackage{mathtools}

\newtheorem{example}{Example}

\newtheorem{theorem}{Theorem}

\newtheorem{remark}{Remark}

\newlength{\noteWidth}
\long\def\notes#1{\ifinner
             {\tiny #1}
             \else
              \marginpar{\parbox[t]{\noteWidth}{\raggedright\tiny #1}}
               \fi}

\def\notes#1{\typeout{#1 !!!}}  

\def\Ebox#1#2{%
\centerline{\includegraphics[scale = #1]{figures/#2}}}

\newcounter{rmnum}

\newenvironment{romannum}{\begin{list}{{\upshape (\roman{rmnum})}}{\usecounter{rmnum}
\setlength{\leftmargin}{6pt}
\setlength{\rightmargin}{4pt}
\setlength{\itemindent}{-1pt}
}}{\end{list}}

\newcounter{anum}

\def\Sec#1{Sec.~\ref{#1}}

\def\IEEEQEDclosed{\mbox{\rule[0pt]{1.3ex}{1.3ex}}}

\def\qed{\nobreak\hfill\IEEEQEDclosed}
\def\UZ{\underline{\mathcal{Z}}}

\def\E{{\sf E}}
\def\P{{\sf P}}

\newcommand{\indicator}[1]{\mathbf{1}_{\left[ {#1} \right] }}

\def\Re{\mathbb{R}}

\def\Inov{I}
\def\v{{\sf K}}
\def\Prob{{\sf P}}

\def\clL{{\cal L}}

\title{\LARGE \bf
Joint Probabilistic Data Association-Feedback Particle Filter \\ for
Multiple Target Tracking Applications}

\author{Tao Yang, Geng Huang, Prashant G. Mehta
\thanks{Financial support from the AFOSR grant FA9550-09-1-0190 and the NSF grant EECS-0925534 is gratefully acknowledged.}
\thanks{T. Yang and P. G. Mehta are with the Coordinated Science Laboratory and
the Department of Mechanical Science and Engineering at the University of
Illinois at Urbana-Champaign (UIUC)
{\tt\scriptsize taoyang1@illinois.edu; mehtapg@illinois.edu}}
\thanks{G. Huang is with the Department of Electrical and Computer
Engineering at UIUC
{\tt\scriptsize huang125@illinois.edu }}%
}

\begin{document}

\maketitle
\thispagestyle{empty}

\begin{abstract}
This paper introduces a novel feedback-control based particle
filter for the solution of the filtering problem with data
association uncertainty.  The particle filter is referred to as
the joint probabilistic data association-feedback particle
filter (JPDA-FPF).  

The JPDA-FPF is based on the feedback particle filter
introduced in our earlier
papers~\cite{YangMehtaMeyn11ACC,YangMehtaMeyn11CDC}. The
remarkable conclusion of our paper is that the JPDA-FPF
algorithm retains the innovation error-based feedback structure
of the feedback particle filter, even with data association
uncertainty in the general nonlinear case. The theoretical
results are illustrated with the aid of two numerical example
problems drawn from multiple target tracking applications.

\end{abstract}

\section{INTRODUCTION}
\label{sec:intro}

Filtering with data association uncertainty is important to a
number of applications, including, air and missile defense
systems, air traffic surveillance, weather surveillance, ground
mapping, geophysical surveys, remote sensing, autonomous
navigation and robotics~\cite{Bar-Shalom_IEEE_CSM,Blackman_book}.  
In each of these applications, there exists data
association uncertainty in the sense that one can not assign
individual measurements to individual targets in an apriori
manner.

Given the large number of applications, algorithms for filtering problems with data
association uncertainty have been extensively studied in the past;
cf.,~\cite{Bar-Shalom_IEEE_CSM,Bar-Shalom_Proc_IEEE} and references therein. A typical
algorithm is comprised of two parts:
\begin{romannum}
\item A {\em filtering} algorithm
  for tracking a single target, and
\item A {\em data association} algorithm for associating
  measurements to targets.
\end{romannum}

Prior to mid-1990s, the primary tool for filtering was a Kalman
filter or one of its extensions, e.g., extended Kalman filter.
The limitations of these tools in applications arise on account
of nonlinearities, not only in dynamic motion of targets (e.g.,
drag forces in ballistic targets) but also in the measurement
models (e.g., range or bearing).  The nonlinearities can lead
to a non-Gaussian multimodal conditional distribution.  For
such cases, Kalman and extended Kalman filters are known to
perform poorly; cf.,~\cite{Ristic_book_2004}. Since the advent
and wide-spread use of particle
filters~\cite{gorsalsmi93,DouFreGor01}, such filters are
becoming increasing relevant to single and multiple target
tracking applications; cf.,~\cite{Ristic_book_2004} and
references therein.

The other part is the data association algorithm. The purpose
of the data association algorithm is to assign measurements to
targets.  The complications arise due to multiple non-target
specific measurements (due to multiple targets in the coverage
area), missing measurements (probability of detection less than
one, e.g., due to target occlusion), false alarms (due to
clutter) and apriori unknown number of targets (that require
track initiation).

The earlier solutions considered assignments in a deterministic
manner: These include the simple but non-robust ``nearest
neighbor'' assignment algorithm and the multiple hypothesis
testing (MHT) algorithm, requiring exhaustive
enumeration~\cite{Reid79analgorithm,Blackman_book}. However,
exhaustive enumeration leads to an NP-hard problem because
number of associations increases exponentially with time.


The complexity issue led to development of probabilistic
approaches: These include the probabilistic MHT or its simpler
``single-scan'' version, the joint probabilistic data
association filter
(JPDAF)~\cite{Bar-Shalom_Proc_IEEE,Bar-Shalom_book_88}.
The central object of interest in these approaches is the
computation (or approximation) of the {\em
measurement-to-target association probability}.  Certain
modeling assumptions are necessary to compute these in a
tractable fashion.  Although the probabilistic algorithms have
reduced computational complexity, they have primarily been
developed for linear settings; cf.,~\cite{Bar-Shalom_IEEE_CSM}.


The development of particle filters has naturally led to
investigations of data association algorithms based on
importance sampling techniques.  This remains an active area of
research; cf.,~\cite{Herman_2011} and references therein. One
early contribution is the multitarget particle filter (MPFT)
in~\cite{Hue00trackingmultiple}.  The data association problem
is approached in the same probabilistic spirit as the basic
JPDAF.  The association probabilities are obtained via the use
of Markov Chain Monte-Carlo (MCMC) techniques; see
also~\cite{Kyriakides_08,Oh_09} for related approaches.

\smallskip

In this paper, we introduce a novel feedback control-based
particle filter algorithm
for solution of the joint filtering-data association
problem.  The proposed algorithm
is based on the feedback particle filter concept introduced by
us in earlier
papers~\cite{YangMehtaMeyn11ACC,YangMehtaMeyn11CDC}. A feedback
particle filter is a controlled system to approximate the
solution of the nonlinear filtering problem. The filter has a
feedback structure similar to the Kalman filter: At each time
$t$, the control is obtained by using a proportional gain
feedback with respect to a certain modified form of the
innovation error.  The filter design amounts to design of the
proportional gain -- the solution is given by the Kalman gain
in the linear Gaussian case.

In the present paper, we extend the feedback particle filter to
problems with data association uncertainty.  We refer to the
resulting algorithm as the joint probabilistic data
association-feedback particle filter (JPDA-FPF).  As the name
suggests, the proposed algorithm represents a generalization of
the Kalman filter-based joint probabilistic data association
(JPDAF) now to the general nonlinear filtering problem.

One remarkable conclusion of our paper is that the JPDA-FPF
retains the innovation error-based feedback structure even for
the nonlinear problem.
The innovation error-based
feedback structure is expected to be useful because of the coupled
nature of filtering and the data association problems.

The theoretical results are illustrated with two numerical
examples. The first example considers the problem of tracking a
single target in the presence of clutter. The second example
considers a multiple target tracking problem. The model problem
scenario is used to illustrate the performance of JPDA-FPF
vis-a-vis possible track coalescence in the presence of data
association uncertainty~\cite{BlomBloem2000}.

The outline of this paper is as follows: The JPDA-FPF algorithm is first described
for single target in the presence of clutter, in
\Sec{sec:PDAFPF}.  The multiple target case follows as a direct
extension, and is discussed in \Sec{sec:JPDAFPF}.  Numerical
examples appear in \Sec{sec:numerics}.

The theory of feedback particle filter can be found in our earlier papers~\cite{YangMehtaMeyn11ACC},\cite{YangMehtaMeyn11CDC}. In the remainder of this paper, we restrict ourselves to the scalar
filtering problem.  The scalar case serves the pedagogical
purpose, and is notationally convenient. The extension to the
multivariable case is straightforward: The feedback particle
filter has the same innovation error-based structure except
that the gain function is replaced by the gain vector.  The
multivariable feedback particle filter is used in the two
numerical examples described in \Sec{sec:numerics}.

\section{Feedback Particle Filter with Data Association Uncertainty}
\label{sec:PDAFPF}


In this section, we describe the probabilistic data association-feedback 
particle filter (PDA-FPF) for the problem of filtering a single target
with multiple measurements. The filter for multiple independent targets is obtained as a
straightforward extension, and briefly described in
\Sec{sec:JPDAFPF}.

\subsection{Problem statement, Assumptions and Notation}


The following notation is adopted:
\begin{romannum}
  \item At time $t$, the target state is denoted by
      $X_t\in\Re$.
  \item At time $t$, the observation vector $
      \underline{Z}_t := ( {Z}_t^1, {Z}_t^2, \hdots,
      {Z}_t^M)^T$, where $M$ is assumed fixed and
      $Z_t^m\in\Re$.
  \item At time $t$, the association random variable
      is denoted as $A_t\in \{0,1,\hdots,M\}$.  It is used
      to associate one measurement to the target: $A_t=m$
      signifies that the $m^{\text{th}}$-measurement
      $Z_t^m$ is 'associated' with the target, and $A_t=0$
      means that the target is not detected at time $t$.
\end{romannum}

The following models are assumed for the three stochastic processes:
\begin{romannum}
\item The state $X_t$ evolves according to a nonlinear
    SDE:
\begin{equation}
\ud X_t = a(X_t)\ud t + \sigma_B \ud B_t,
\label{eqn:Signal_Process_Target}
\end{equation}
where $a(\cdot)$ is a $C^1$ function and $\{B_t\}$ is a standard Wiener process.
\item The association random process $A_t$ evolves as
    a jump Markov process in continuous-time:
\begin{equation}
{\sf P}(A_{t+\delta}=m'|A_{t}=m) = \frac{c}{M} \delta +
o(\delta),\quad m'\ne m.
\end{equation}
The initial distribution ${\sf P}([A_0=m])=\frac{1}{M+1}$. $A_t$ and $X_t$ are assumed to be mutually
    independent.
\item At time $t$, the observation model is given by,
\begin{equation}
\ud Z_t^{m} = 1_{[A_t=m]} h(X_t) \ud t + \sigma_W \ud W_t^m,\label{eqn:observ_model_seperate}
\end{equation}
for $m\in\{1,\hdots,M\}$, where $h(\cdot)$ is $C^1$ function and $\{W_t^m\}$ are mutually
independent standard Wiener processes and
\begin{equation*}
1_{[A_t=m]} :=\begin{dcases*}
        1  &  when $A_t=m$ \\
        0 & otherwise.
        \end{dcases*}
\end{equation*}

\end{romannum}

The problem is to obtain the posterior distribution of
${X}_t$ given the history of observations $\mathcal{\underline{Z}}_t :=
\sigma(\underline{Z}_s:s\le t)$. 



The methodology comprises of the following two parts:
\begin{romannum}
\item Evaluation of association probability, and
\item Integration of association probability in the
    feedback particle filter framework.
\end{romannum}

\subsection{Association Probability for a Single Target}

The association probability is defined as the probability of
the association $[A_t=m]$ conditioned on
$\mathcal{\underline{Z}}_t$:
\begin{equation}
\beta_t^m \triangleq \Prob ([A_t=m]|\mathcal{\underline{Z}}_t),\quad m = 0,1,...,M. \label{eqn:def_betaj}
\end{equation}
Since the events are mutually exclusive and exhaustive,
$\sum_{m=0}^{M} \beta_t^m = 1$.

For the single-target-multiple-observation model described
above, the filter for computing association probability is
derived in Appendix~\ref{apdx:association_filter_pda}.
It is of the following form: For $m\in\{1,...,M\}$,
\begin{align}
\ud \beta_t^m & = \frac{c}{M} \left[ 1 - (M+1) \beta_t^m \right] \ud t \nonumber\\
& + \frac{1}{\sigma_W^2}\beta_t^m \hat{h}_t \sum_{j=1}^M \beta_t^j \left[(\ud Z_t^m - \beta_t^m \hat{h}_t\ud t)- (\ud Z_t^j - \beta_t^j \hat{h}_t\ud t)\right]\nonumber\\
& + \frac{1}{\sigma_W^2}\beta_t^m  (\widehat{h_t^2} - \hat{h}_t^2) \sum_{j=1}^M \beta_t^j (\beta_t^j - \beta_t^m)\ud t,
\label{eqn:filter_for_beta_nonlinear}
\end{align}
where $\hat{h}_t =  {\sf E} [h(X_t)|\mathcal{Z}_t]$ and
$\widehat{h^2_t} = {\sf E} [h^2(X_t)|\mathcal{Z}_t]$. These are
approximated by using the feedback particle filter described in
the following section.

In practice, one may also wish to consider approaches to reduce
filter complexity, e.g., by assigning gating regions for the
measurements;  cf., Sec.~4.2.3 in~\cite{Bar-Shalom_book_88}.

\begin{remark}
In the following, we integrate association probability with the
feedback particle filter, which is used to approximate
evolution of the posterior.  Separate algorithms for data
association and posterior are motivated in part by the
classical JPDA filtering
literature~\cite{Bar-Shalom_book_88,Bar-Shalom_IEEE_CSM,Bar-Shalom_Proc_IEEE}.
A separate treatment is also useful while considering multiple
target tracking problems.  For such problems, one can extend
algorithms for data association in a straightforward manner,
while the algorithm for posterior remains as before. This is
illustrated with the aid of two-target-two-observation example
in Sec~\ref{sec:JPDAFPF}.
\end{remark}

\begin{remark}
The association probability
filter~\eqref{eqn:filter_for_beta_nonlinear} can also be
derived by considering a continuous-time limit starting from the
continuous-discrete time filter in literature~\cite{Bar-Shalom_book_88}.
This proof appears
in Appendix~\ref{apdx:discretized_assoc_filter}.  The alternate proof
is included for the following reasons:
\begin{romannum}
\item The proof shows that the
    filter~\eqref{eqn:filter_for_beta_nonlinear}
  is in fact the continuous-time nonlinear counterpart of the
  algorithm that is used to
  obtain association probability in the classical JPDAF filter.  This is
  important because some of the modeling assumptions (e.g., modeling
  of clutter, or of association $A_t$ via a jump Markov process) here may appear
  to be different from those considered in the classical literature.
\item The proof method suggests alternate discrete-time algorithms for evaluating
  association probabilities in simulations and experiments, where
  observations are made at discrete sampling times.
\end{romannum}
\end{remark}

\subsection{Probabilistic Data Association-Feedback Particle Filter}

Following the feedback particle filter methodology, the model for the particle filter is
given by,
\begin{equation}
\ud X^i_t = a(X^i_t) \ud t + \sigma_B \ud B^i_t  + \ud U^i_t , \label{eqn:particle_model}
\end{equation}
where $X^i_t \in \Re$ is the state for the $i^{\text{th}}$
particle at time $t$, $ U^i_t$ is its control input, and
$\{B^i_t\}$ are mutually independent standard Wiener processes.
We assume the initial conditions $\{X^i_0\}_{i=1}^N$  are
i.i.d., independent of $\{B^i_t\}$, and drawn from the initial
distribution $p^*(x,0)$ of $X_0$.  Both  $\{B^i_t\}$ and
$\{X^i_0\}$ are also assumed to be independent of $X_t,Z_t$.
Certain additional assumptions are made regarding admissible
forms of control input (see~\cite{YangMehtaMeyn11CDC}).

Recall that there are two types of conditional distributions of
interest in our analysis:
\begin{romannum}
\item $p(x,t)$:  Defines the conditional dist.\ of
    $X^i_t$ given $\UZ_t$.
\item $p^*(x,t)$: Defines the conditional dist.\ of
    $X_t$ given $\UZ_t$.
\end{romannum}

The control problem is to choose the control input $U^i_t$ so
that $p$ approximates $p^*$, and consequently empirical
distribution of the particles approximates $p^*$ for large
number of particles.

The evolution of  $p^*(x,t)$ is described by modified form of the Kushner-Stratonovich (K-S) equation:
\begin{equation}
\ud p^\ast = \clL^\dagger p^\ast \ud t +
\frac{1}{\sigma_W^2}
\sum_{m=1}^{M} \beta_t^m ( h-\hat{h}_t )(\ud Z_t^m - \hat{h}_t \ud t)p^\ast. \label{eqn:PDA-KS}
\end{equation}
where $ \hat{h}_t = \int h(x) p^*(x,t) \ud x$, and $
\clL^\dagger$ is the Kolmogorov forward operator. The proof
appears in Appendix~\ref{apdx:consistency}.

The main result of this section is to describe an explicit
formula for the optimal control input, and demonstrate that
under general conditions we obtain an exact match:  $p=p^*$
under optimal control. The optimally controlled dynamics of the
$i^{\text{th}}$ particle have the following form,
\begin{align}
\ud X_t^i = a(X^i_t)\ud t + &\sigma_B \ud B_t^i + \sum_{m=1}^{M} \beta_t^m \v(X_t^i,t) \ud \Inov^{i,m}_t \nonumber\\
&+ \frac{1}{2}\sigma_W^2 \sum_{m=1}^{M} (\beta_t^m)^2 \v(X_t^i,t) \v'(X_t^i,t) \ud t,\label{eqn:PDA-FPF}
\end{align}
where $\Inov^{i,m}_t$ is a modified form of the
$\emph{innovation process}$,
\begin{equation}
\ud \Inov^{i,m}_t := \ud Z_t^m - [\frac{\beta_t^m}{2} h(X_t^i) + (1-\frac{\beta_t^m}{2})\hat{h}_t]\ud t,\label{eqn:pda-inov}
\end{equation}
where $\hat{h}_t := {\sf E} [h(X_t)|\mathcal{Z}_t] = \int h(x)
p(x,t) \ud x$.  The gain function $\v$ is the solution of a certain
EL-BVP:
\begin{equation}
-\frac{\partial}{\partial x}\left( \frac{1}{p(x,t)} \frac{\partial}{\partial x} \{ p(x,t) \v(x,t) \}\right) = \frac{1}{\sigma_W^2}h'(x),
\label{eqn:ELBVP_single}
\end{equation}

The evolution of $p(x,t)$ is easily obtained as the forward
Kolmogorov operator: See Appendix~\ref{apdx:consistency} for
the equations.


The following theorem shows that the two evolution equations
for $p$ and $p^\ast$ are identical. The proof appears in
Appendix~\ref{apdx:consistency}.

\begin{theorem}\label{thm:kushner}
Consider the two evolutions for $p$ and $p^\ast$, defined
according to the Kolmogorov forward equation and modified K-S
equation~\eqref{eqn:PDA-KS}, respectively. Suppose that the
gain function $\v(X,t)$ is obtained according
to~\eqref{eqn:ELBVP_single}. Then provided $p(x,0) =
p^\ast(x,0)$, we have for all $t \leq 0$, $p(x,t) =
p^\ast(x,t)$. \qed
\end{theorem}

\begin{example} Consider the
single target, single measurement case where the measurement
may be due to clutter (false alarm). Let $\beta_t$ denote the
measurement-to-target association probability at time $t$.

For this case, the feedback particle filter is given by the
controlled system which is a special case
of~\eqref{eqn:PDA-FPF},
\begin{align}
\ud X^i_t &= a(X^i_t) \ud t + \sigma_B \ud B^i_t \nonumber\\
          &\quad + \beta_t \v(X^i_t,t) \ud \Inov^i_t + \frac{1}{2} \beta_t^2 \sigma_W^2 \v(X^i_t,t) \v'(X^i_t,t) \ud t,
\label{eqn:particle_filter_clutter}
\end{align}
where the innovation error $\Inov^i_t$ is given by,
\begin{equation}
\ud \Inov^i_t := \ud Z_t - \left(\frac{\beta_t}{2} h(X^i_t) +
  (1-\frac{\beta_t}{2}) \hat{h}_t \right) \ud t.
\label{e:in_clutter}
\end{equation}

\noindent For the two extreme values of $\beta_t$, the filter
reduces to the known form:
\begin{romannum}
\item If $\beta_t=1$, the measurement is associated
    with the target with probability $1$. In this case, the
    filter is the same as
    FPF presented
    in~\cite{YangMehtaMeyn11CDC}.
\item If $\beta_t=0$, the measurement carries no
    information and the control input $\ud U^i_t =0$.
\end{romannum}
For $\beta_t\in(0,1)$, the control is more interesting. The
remarkable fact is that the innovation error-based feedback
control structure is preserved.  The association probability
serves to modify the formulae for the gain function and the
innovation error:
\begin{romannum}
\item The gain function is effectively reduced to
    $\beta_t \v(X^i_t,t)$.  That is, the control gets less
    agressive in the presence of possible false alarms due
    to clutter.
\item The innovation error is given by a more
    general formula~\eqref{e:in_clutter}.  The optimal
    prediction of the $i^{\text{th}}$-particle is now a
    weighted average of  $h(X^i_t)$ and the population
    prediction
$\hat{h}_t \approx \frac{1}{N} \sum_{j=1}^N h(X^j_t)$.
  Effectively, in the presence of possible false alarms, a
  particle gives more weight to the population in computing
  its innovation error.
\end{romannum}
\end{example}

\subsection{Example: Linear Case}
We provide here a special case of PDA-FPF for the single target
tracking problem described by a linear model:
\begin{subequations}
\begin{align}
\ud X_t  &= \alpha \;X_t\ud t + \sigma_B \ud B_t,\label{eqn:dyn_lin}\\
\ud Z_t &= \gamma \; X_t \ud t+\sigma_W \ud W_t,\label{eqn:obs_lin}
\end{align}
\end{subequations}
where $\alpha,\gamma$ are real numbers. 

The PDA-FPF is described by~\eqref{eqn:PDA-FPF}--\eqref{eqn:ELBVP_single}. If we assume $p(x,t)$ to be Gaussian at each time with mean $\mu_t$ and variance $\Sigma_t$, i.e., $p(x,t)= \frac{1}{\sqrt{2 \pi \Sigma_t}}
\exp(-\frac{(x-\mu_t)^2}{2\Sigma_t})$, then by direct substitution in~\eqref{eqn:ELBVP_single} we obtain the gain function:
\begin{equation}
\v(x,t) = \frac{\gamma\Sigma_t}{\sigma_W^2}.
\label{eqn:linsol_v}
\end{equation}


The PDA-FPF is then given by,
\begin{align}
\ud X^i_t &= \alpha \; X^i_t \ud t + \sigma_B \ud B^i_t\nonumber\\
&+ \frac{\gamma \Sigma_t}{\sigma^2_W} \sum_{m=1}^{M} \beta_t^m \left[ \ud Z_t^m - \gamma (\frac{\beta_t^m}{2} X^i_t + (1-\frac{\beta_t^m}{2}) \mu_t) \ud t \right].
\label{eqn:pda_fpf_lin}
\end{align}

%

The filter for association probability $\beta_t^m$ is as follows:
\begin{align}
\ud &\beta_t^m  = \frac{c}{M} \left[ 1 - (M+1) \beta_t^m \right] \ud t \nonumber\\
& + \frac{1}{\sigma_W^2}\beta_t^m \gamma \mu_t \sum_{j=1}^M \beta_t^j \left[(\ud Z_t^m - \beta_t^m \gamma \mu_t \ud t)- (\ud Z_t^j - \beta_t^j \gamma \mu_t \ud t)\right]\nonumber\\
& + \frac{1}{\sigma_W^2}\beta_t^m \gamma^2 \Sigma_t \sum_{j=1}^M \beta_t^j (\beta_t^j - \beta_t^m)\ud t,
\label{eqn:filter_for_beta_linear}
\end{align}

In practice $\{\mu_t, \Sigma_t\}$ in \eqref{eqn:pda_fpf_lin}-\eqref{eqn:filter_for_beta_linear} are approximated as sample
means and sample covariances using
$\{X^i_t\}_{i=1}^N$.
\begin{equation}
\begin{aligned}
\mu_t & \approx \mu_t^{(N)} := \frac{1}{N} \sum_{i=1}^N X^i_t,\\
\Sigma_t & \approx \Sigma_t^{(N)} := \frac{1}{N-1} \sum_{i=1}^N (X^i_t - \mu_t^{(N)})^2.
\end{aligned}
\label{e:mut_sigmat_approx}
\end{equation}


\section{Multiple Target Tracking using \\Feedback Particle Filter}
\label{sec:JPDAFPF}

In this section, we extend the PDA-FPF to multiple target
tracking problems. Specifically, a two-target two-observation
problem is used to illustrate JPDA-FPF. The extension to the
more general case is straightforward.

\subsection{Problem statement}

The following notation is adopted:
\begin{romannum}
\item At time $t$, the target state is denoted as
    $\underline{X}_t:=(X_t^1,X_t^2)^T$, where
    $X_t^n\in\Re$ for $n \in \{1,2\}$.
\item   At time $t$, the observation vector $
    \underline{Z}_t := ({Z}_t^1, {Z}_t^2)^T$, where
    $Z_t^m\in\Re$ for $m\in\{1,2\}$.
\item At time $t$, the association random variable
    is denoted as $A_t\in \{1,2\}$.  It is used to
    associate measurements to targets in a joint manner:
    $A_t=1$ signifies that $Z_t^1$ is associated with
    target $1$ and $Z_t^2$ with target $2$. Similarly
    $A_t=2$ accounts for the complementary case.
\end{romannum}

The following models are assumed for the three stochastic
processes:
\begin{romannum}
\item Each element of the state vector
    $\underline{X}_t$ evolves according to a
    one-dimensional nonlinear SDE:
\begin{equation}
\ud X_t^n = a(X_t^n)\ud t + \sigma_B^n \ud B_t^n,\quad n\in\{1,2\}
\label{eqn:Signal_Process_Two_Target}
\end{equation}
where $\{B_t^1\}$,$\{B_t^2\}$ are mutually independent standard
Wiener processes.
\item The association random process $A_t$ evolves as
    a jump Markov process in continuous-time:
\begin{equation}
{\sf P}(A_{t+\delta}=m'|A_{t}=m) = c \delta +
o(\delta),\quad m'\ne m.
\end{equation}
The initial distribution ${\sf P}([A_0=m])=\frac{1}{2}$. $A_t$ and $\underline{X}_t$ are assumed to be
    mutually independent.
\item At time $t$, the observation model is given by,
\begin{equation}
\left[
\begin{array}{ccc}
\ud Z_t^1 \\
\ud Z_t^2
\end{array} \right] = \Psi (A_t) \left[
                                    \begin{array}{ccc}
                                    h(X_t^1)\\
                                    h(X_t^2)
                                    \end{array}\right] \ud t +\sigma_W\left[
                                                         \begin{array}{ccc}
                                                         \ud W_t^1 \\
                                                         \ud W_t^2
                                                         \end{array}\right],\label{eqn:Two_Target_Observ_Process}
\end{equation}
where $\{W_t^1\},\{W_t^2\}$ are mutually independent standard
Wiener processes and $\Psi(A_t)$ is a function which maps $A_t$
to a permutation matrix:
\begin{equation}
\Psi(1) =
\begin{bmatrix}
        1 & 0 \\
        0 & 1
       \end{bmatrix},\quad \Psi(2) = \begin{bmatrix}
                                                0 & 1\\
                                                1 & 0
                                                \end{bmatrix}.\label{eqn:permuation}
\end{equation}
\end{romannum}

\subsection{Joint Probabilistic Data Association for Two Target}

The joint association probability is defined as the probability
of the joint association $[A_t = m]$ conditioned on $\UZ_t$:
\begin{equation}
\pi_t^m \triangleq \P \left([A_t = m]|\UZ_t\right),\quad m=1,2. \label{eqn:jpda_def}
\end{equation}

The filter for joint association probability $\pi_t$ is a
straightforward extension
of~\eqref{eqn:filter_for_beta_nonlinear}. The proof appears in
Appendix~\ref{apdx:association_filter_jpda}. It is of the
following form:
\begin{align}
&\ud \pi_t^1 = -c (\pi_t^1 - \pi_t^2)\ud t + \frac{1}{\sigma_W^2}\pi_t^1 \pi_t^2 (\hat{h}^1_t-\hat{h}_t^2)(\ud \mu_t^1 - \ud \mu_t^2)\nonumber\\
& -\frac{1}{\sigma_W^2}\pi_t^1\pi_t^2 (\pi_t^1 - \pi_t^2)[\widehat{(h_t^1)^2} - (\hat{h}_t^1)^2 + \widehat{(h_t^2)^2} - ((\hat{h}_t^2)^2 ] \ud t,
\label{eqn:jpda_dpi_1}
\end{align}
where $\ud \mu_t^m = \ud Z_t^m - (\pi_t^1 - \pi_t^2)
\hat{h}_t^m \ud t$, $\hat{h}_t^n:= \E [h(X_t^n)|\UZ_t]$ and
$\widehat{(h_t^n)^2}:= \E[h^2(X_t^n)|\UZ_t]$. Since the joint
events are mutually exclusive and exhaustive, we have
$\sum_{m=1}^2 \pi_t^m = 1$. Using this, we have $\pi_t^2 = 1-
\pi_t^1$ and $\ud \pi_t^2 = - \ud \pi_t^1$.

\subsection{Joint Prob. Data Association-Feedback Particle Filter}
The joint association probabilities $\pi_t^1,\pi_t^2$ are used
to obtain marginal association probability for individual
target. For example, for target $1$: $\beta_t^1 = \pi_t^1$,
$\beta_t^2 = \pi_t^2$. Once the association probabilities are
known, the feedback particle filter for each target is of the
form~\eqref{eqn:PDA-FPF}.

\section{Numerics}
\label{sec:numerics}

In this section, we discuss results of two numerical examples.
Even though the theory was described for real-valued state and
observation processes, the numerical examples consider more
realistic multivariable models.

\subsection{Single Target Tracking in Clutter}
\label{subsec:stt}

We first consider a single target tracking problem where the
target dynamics evolve according to a white-noise acceleration
model:
\begin{align}
\ud X_t &= FX_t\ud t + \sigma_B \ud B_t,\label{eqn:stt_dyn}\\
\ud Z_t &= HX_t\ud t + \sigma_W \ud W_t,\label{eqn:stt_obs}
\end{align}
where $X_t$ denotes the state vector comprising of position and
velocity coordinates at time $t$, $Z_t$ is the observation
process, $\{B_t\},\{W_t\}$ are mutually independent standard
Wiener processes. The two matrices are given by:
\begin{equation}
F=\begin{bmatrix}
                    0 &  1 \\
                    0 &  0
  \end{bmatrix},
\quad H = \begin{bmatrix}
                    1 &  0
    \end{bmatrix}.
\end{equation}

In the simulation results described next, we use the following
parameter values: $\sigma_B = [0 ;1]$,
$\sigma_W = 0.06$ and initial condition $X_0 =
[0;6]$.  The total simulation time is
$T=1$ and time step $\ud t = 0.01$. At each
discrete-time step, we assume $M=4$ measurements, one due to
target and other three due to clutter. The associations are not
apriori known.

\begin{figure}
\Ebox{.28}{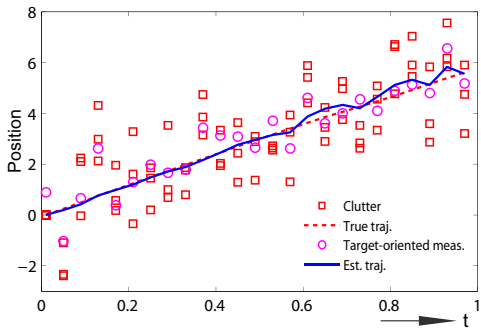}
\caption{Simulation results of single target tracking in clutter using PDA-FPF: Comparison of estimated mean with the true trajectory.}
\label{fig_stt}
\end{figure}

Figure~\ref{fig_stt} depicts the result of a single simulation:
True target trajectory is depicted as a dashed line. At each
discrete time step, target-oriented measurements are depicted
as circles while clutter measurements are depicted as squares.
The estimated mean trajectory is depicted as a solid line. It
is obtained using the PDA-FPF algorithm described in
Sec~\ref{sec:PDAFPF}-D. For the filter simulation, we use
$N=1000$ particles.


\subsection{Track Coalescence Avoidance using JPDA-FPF}
\label{subsec:col_avoid}

Track coalescence is a common problem in multiple tracking
applications. Track coalescence can occur when two closely
spaced targets move with approximately the same velocity over a
time period~\cite{Blackman_book}. With standard
implementations of JPDAF and SIR particle filter algorithms,
the target tracks tends to coalesce even after the targets have
moved apart~\cite{BlomBloem2000}. In the following example, we
describe simulation results for JPDA-FPF for a model problem
scenario taken from~\cite{BlomBloem2006}.

We consider two targets. For each target, the dynamics are
described by a white-noise acceleration model as in the
preceding example. For $n=1,2$:
\begin{equation}
\ud X_t^n = F X_t^n \ud t + \sigma_B \ud B_t, \label{eqn:col_dyn}
\end{equation}
where the state $X_t^n$ comprises of target position and
velocity.

We assume two observations are given by
\begin{equation}
\left[
\begin{array}{ccc}
\ud Z_t^1 \\
\ud Z_t^2
\end{array} \right] = \Psi (A_t) \left[\begin{array}{ccc}
                                        H X_t^1\\
                                        H X_t^2
                                        \end{array}\right] \ud t +\sigma_W\left[
                                                         \begin{array}{ccc}
                                                         \ud W_t^1 \\
                                                         \ud W_t^2
                                                         \end{array}\right],\label{eqn:col_obs}
\end{equation}
where $A_t$ is the association random variable, $\Psi(A_t)$ is
the permutation matrix as defined in~\eqref{eqn:permuation}. $\{B_t\},\{W_t^1\},\{W_t^2\}$ are mutually independent
standard Winer processes.

In the simulation results described next, we use the following
parameter values: $\sigma_B = [0;2]$,
$\sigma_W = [0.005;0.005]$ and initial
condition $\underline{X}_0=[1;
-3.5;-1;3.5]$. The total simulation
time is $T=1s$ and time step $\ud t=0.001s$. The prior
association probability $(\pi_t^1,\pi_t^2)$ is assumed to be
$(1/2,1/2)$.

\begin{figure*}
\centering \hspace{-0.2cm}
\includegraphics[scale = 0.25]{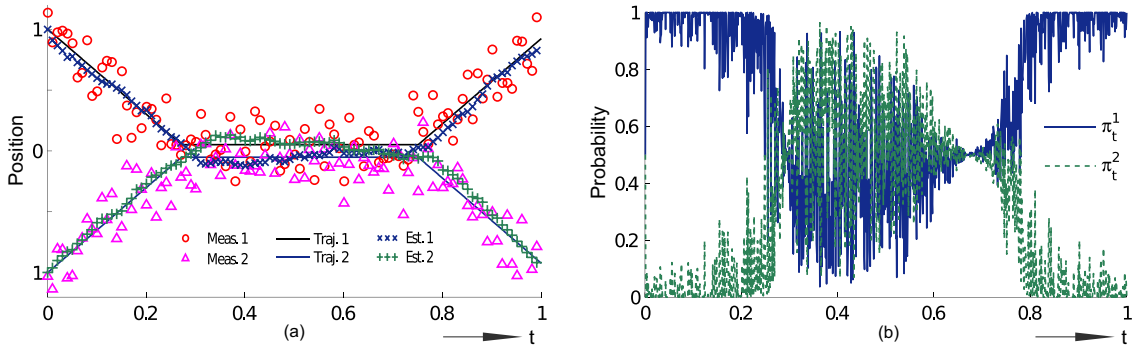}

 \caption{Track coalescence
avoiding using JPDA-FPF: (a) Comparison of estimated mean
JPDA-FPF with true trajectories. (b) Plot of data association
probability.} \vspace{0.25cm} \label{fig_col_track}
\end{figure*}

Figure~\ref{fig_col_track}(a) depicts the results of a single
simulation:  The estimated mean trajectories are obtained using
the JPDA-FPF described in Sec~\ref{sec:JPDAFPF}.
Figure~\ref{fig_col_track}(b) depicts the evolution of
association probability $(\pi_t^1,\pi_t^2)$ during the same
simulation run. For the filter simulation, we use $N=1000$
particles. To obtain the association probabilities, we use an
adaptive time stepping scheme for numerical integration of
association probability filter~\eqref{eqn:jpda_dpi_1}.



%

\section*{Acknowledgement}
We are grateful to Prof. Sean Meyn for useful discussions on the work described in this paper.


\appendix

\subsection{Association Probability Filter for $\beta_t^m$}
\label{apdx:association_filter_pda}

Rewrite~\eqref{eqn:observ_model_seperate} in the vector form:
\begin{equation}
\ud \underline{Z}_t = \underline{\chi}(A_t) h(X_t)\ud t+ \sigma_W \ud \underline{W}_t,\label{eqn:observ_model_vector}
\end{equation}
where $\underline{\chi}(A_t) := (\chi_t^1,...,\chi_t^M)^T$,
$\chi_t^m = \indicator{A_t=m}$, and
$\underline{W}_t:=(W_t^1,...,W_t^M)^T$. The transition
intensity matrix for the jump Markov process $A_t$ is denoted
as $\Lambda$ with
\begin{equation}
\Lambda_{mm'} = \begin{dcases*}
        -c  &  if $m=m'$\\
        c/M &  if $m\neq m'$.
        \end{dcases*}
\end{equation}
Denote $\mathcal{X}_t:= \sigma(X_s:s \leq t)$, and
$\mathcal{C}_t := \mathcal{X}_t \vee \mathcal{Z}_t$. The
derivation is based on the property of the conditional
expectation:
\begin{equation}
\E[\underline{\varphi}(A_t)|\mathcal{Z}_t] = \E[\E[\underline{\varphi}(A_t)|\mathcal{C}_t]|\mathcal{Z}_t].\nonumber
\end{equation}
The SDE for evolution of
$\E[\underline{\varphi}(A_t)|\mathcal{C}_t]$ is described by
the standard Wonham filter:
\begin{align}
&\E[\underline{\varphi}(A_t)|\mathcal{C}_t] = \E[\underline{\varphi}(A_0)] + \int_0^t \E[\Lambda \underline{\varphi}(A_s)|\mathcal{C}_s]\ud s\nonumber\\
&+ \sum_{m=1}^M \int_0^t \E[(D_s^m - h(X_s)\beta_s^m I) \underline{\varphi}(A_s)|\mathcal{C}_s] (\ud Z_s^m - h(X_s)\beta_s^m \ud s),
\label{eqn:proof_pda_wonham}
\end{align}
where $I$ is the identity matrix, $D_t^m$ is the diagonal
matrix where the only non-zero entry is $(D_t^m)_{mm} =
h(X_t)$.

Taking $\E[\cdot|\mathcal{Z}_t]$
of~\eqref{eqn:proof_pda_wonham} gives the desired
result~\eqref{eqn:filter_for_beta_nonlinear}.

\subsection{Consistency proof of $p$ and $p^\ast$}
\label{apdx:consistency}

\noindent {\bf Evolution equation for $p^\ast$:} Recall
$\mathcal{\underline{Z}}_t := \sigma(\underline{Z}_s:s\le t)$,
$\mathcal{A}_t:=\sigma(A_s:s\le t)$. We denote $\mathcal{C}_t
:= \mathcal{A}_t \vee \mathcal{\underline{Z}}_t$.  The
derivation is based on the property of the conditional
expectation:
\begin{equation*}
{\sf E} [ \varphi(X_t)|\mathcal{\underline{Z}}_t] = {\sf E} \left[
  {\sf E} [ \varphi(X_t)|\mathcal{{C}}_t] |
  \mathcal{\underline{Z}}_t \right].
\end{equation*}
The sde for evolution of  ${\sf E} [
\varphi(X_t)|\mathcal{{C}}_t]$ is described by the standard
nonlinear filter with innovation error, $\sum_{m=1}^M \chi_t^m (\ud
Z_t^m - \hat{h}_t \ud t)$, where $\chi_t^m = \indicator{A_t=m}$:
\begin{align*}
{\sf E} [\varphi(X_t)|\mathcal{C}_t] &= {\sf E} [
\varphi(X_0)] +  \int_{0}^t {\sf E} [
\clL \varphi(X_s)|\mathcal{C}_s] \ud s \\
 &+ \sum_{m=1}^M  \int_{0}^t {\sf E} [(
h-\hat{h}_s)\varphi(X_s) |\mathcal{C}_s] \chi_s^m (\ud
Z_s^m - \hat{h}_s \ud s),
\end{align*}
where $\clL$ denotes the Kolmogorov's backward operator for the diffusion~\eqref {eqn:Signal_Process_Target} (the adjoint of  $\clL^\dagger$).

Taking $  {\sf E} [\cdot|\mathcal{\underline{Z}}_t]$ gives the desired
result because ${\sf E} [\chi_s^m |\mathcal{\underline{Z}}_s] = {\sf
  P}([A_s = m]|\mathcal{\underline{Z}}_s) = \beta_s^m$.

\smallskip
\noindent {\bf Evolution equation for $p$:}
We express the FPF~\eqref{eqn:PDA-FPF} as:
\[
\ud X_t^i = a(X^i_t)\ud t + \sigma_B \ud B_t^i +  \v(X_t^i,t) \sum_{m=1}^{M}
\beta_t^m \ud Z^m_t + u(X_t^i,t) \ud t,
\]
where
\begin{align}
u(x,t):=& - \sum_{m=1}^{M} \beta_t^m \left[ \frac{\beta_t^m }{2} h +
  (1 - \frac{\beta_t^m }{2}) \hat{h}_t \right] \v(x,t)  \nonumber \\ &+
\frac{\sigma_W^2}{2} \sum_{j=1}^{M} (\beta_t^m)^2 \v \v'(x,t). \label{eqn:PDA_u_def}
\end{align}
The evolution equation for $p$ now follows:
\begin{align}
\ud p = \clL^\dagger p \ud t & -
\frac{\partial}{\partial x}( u p ) \ud t
+  \frac{\sigma_W^2}{2} \sum_{j=1}^{M}
(\beta_t^m)^2\frac{\partial^2}{\partial x^2}\left( p \v^2 \right) \ud t \nonumber\\
& - \frac{\partial}{\partial x}\left( \v
  p \right) \sum_{m=1}^{M} \beta_t^m \ud Z_t^m.
\label{eqn:mod_PDA_FPK}
\end{align}

\smallskip

\noindent {\bf Proof of consistency.} The proof follows closely
the consistency proof for the feedback particle filter (see
Appendix~C in~\cite{YangMehtaMeyn11CDC}).  If $\v$ solves the
E-L BVP then
\begin{equation}
-  \frac{\partial}{\partial x} (\v p) = \frac{1}{\sigma_W^2} (h-\hat{h}_t) p\label{eqn:pv}
\end{equation}
On multiplying both sides
of~\eqref{eqn:PDA_u_def} by $-p$ and simplifying (by using~\eqref{eqn:pv}),
we obtain
\[
-up = -\frac{\sigma_W^2}{2} \frac{\partial}{\partial x}(p \v^2)
\sum_{m=1}^M (\beta_t^m)^2 + \v p \sum_{m=1}^M \beta_t^m \hat{h}_t
\]
Differentiate now both sides with respect to $x$ and use~\eqref{eqn:pv} once again to arrive at
\begin{equation}
\frac{\sigma_W^2}{2} \sum_{j=1}^{M}
(\beta_t^m)^2 \frac{\partial^2}{\partial x^2}\left( p \v^2 \right) - \frac{\partial}{\partial x}( u p )
= \frac{-1}{\sigma_W^2} (h-\hat{h}_t) p \sum_{m=1}^M \beta_t^m \hat{h}_t. \label{eqn:ns2}
\end{equation}
Using~\eqref{eqn:pv} and~\eqref{eqn:ns2}
in the forward equation~\eqref{eqn:mod_PDA_FPK}, we obtain:
\begin{equation}
\ud p = \clL^\dagger p \ud t + \frac{1}{\sigma_W^2}\sum_{j=1}^{M} \beta_t^m (h-\hat{h}_t)(\ud Z_t^m- \hat{h}_t\ud t)p\label{eqn:PDA_p_FPK}
\end{equation}
This is precisely the SDE~\eqref{eqn:PDA-KS}, as desired.

\subsection{Association Probability filter for $\pi_t^m$}
\label{apdx:association_filter_jpda}

The derivation follows closely the derivation in
Appendix~\ref{apdx:association_filter_pda}. Note that the
observation model is described
by~\eqref{eqn:Two_Target_Observ_Process}. Denote
$\underline{h}(\underline{X}_t) := (h(X_t^1),h(X_t^2))^T$,
$\widehat{\Psi}_t := \sum_{m=1}^2 \pi_t^m \Psi(m)$ and
$\underline{\phi}_t = (\phi_t^1,\phi_t^2)^T:=\widehat{\Psi}_t
\underline{h}(\underline{X}_t)$. The Wonham filter is given by:
\begin{align}
&\E[\underline{\varphi}(A_t)|\mathcal{C}_t] = \E[\underline{\varphi}(A_0)] + \int_0^t \E[\Lambda \underline{\varphi}(A_s)|\mathcal{C}_s]\ud s\nonumber\\
&+ \sum_{m=1}^2 \int_0^t \E[(D_s^m - \phi_s^m I) \underline{\varphi}(A_s)|\mathcal{C}_s] (\ud Z_s^m - \phi_s^m \ud s),
\label{eqn:proof_jpda_wonham}
\end{align}
where $D_t^m$ is a $2\times 2$ diagonal matrix where
$(D_t^m)_{ii}$ is the $m^{\text{th}}$ entry of the vector
$\Psi(i) \underline{h}(\underline{X}_t)$.

Taking $\E[\cdot|\mathcal{Z}_t]$
of~\eqref{eqn:proof_jpda_wonham} gives the desired result.

\subsection{Alternate Derivation of~\eqref{eqn:filter_for_beta_nonlinear}}
\label{apdx:discretized_assoc_filter}

The aim of this section is to derive, formally, the update part
of the continuous time
filter~\eqref{eqn:filter_for_beta_nonlinear} by taking a
continuous time limit of the discrete-time algorithm for
evaluation of association probability. The procedure for taking
the limit is similar to Sec~$6.8$ in~\cite{jazwinski70} for
derivation of the K-S equation.

At time $t$, we have $M$ measurements $\ud \underline{Z}_t =
(\ud Z_t^1,\ud Z_t^2,...,\ud Z_t^M)^T$, only one of which
originates from the target. The discrete-time filter for
association probability is obtained by using Bayes' rule (see
~\cite{Bar-Shalom_book_88}):
\begin{equation}
\P([A_t = m]|\UZ_t,\ud \underline{Z}_t) = \frac{\P(\ud \underline{Z}_t|[A_t = m])\P([A_t = m]|\UZ_t)}{\sum_{j=1}^M\P(\ud \underline{Z}_t|[A_t = j])\P([A_t = j]|\UZ_t)}.\label{eqn:discret_bayes}
\end{equation}

In evaluation of the association probability, one typically
assumes a clutter model whereby the independent measurements
are uniformly and independently distributed in the coverage
area $V$
(\cite{Bar-Shalom_book_88},\cite{Bar-Shalom_IEEE_CSM}). We then
have:
\begin{align}
\P(\ud \underline{Z}_t|[A_t =m]) &=V^{1-M} \P(\ud Z_t^m|[A_t=m]) \nonumber\\
                          &= V^{1-M} L(\ud Z_t^m).\label{eqn:Z_t_m}
\end{align}
where $L(\ud Z_t^m) = \frac{1}{\sqrt{2\pi\sigma_W^2\ud
t}}\int_\mathbb{R} \exp\left[-\frac{(\ud Z_t^m - h(x)\ud
t)^2}{2 \sigma_W^2 \ud t}\right]p(x,t)\ud x$.

Now, denote $\beta_t^m = \P([A_t=m]|\UZ_t)$, the increment in
the measurement update step (see Sec~$6.8$
in~\cite{jazwinski70}) is given by
\begin{equation}
\ud \beta_t^m:= \P([A_t=m]|\UZ_t,\ud Z_t) - \P([A_t=m]|\UZ_t).\label{eqn:def_dbeta}
\end{equation}
Using~\eqref{eqn:discret_bayes} and~\eqref{eqn:def_dbeta}, we
have:
\begin{equation}
\ud \beta_t^m = E^m(\ud t,\ud \underline{Z}_t)\beta_t^m  - \beta_t^m ,\label{eqn:cal_dbeta}
\end{equation}
where
\begin{equation}
E^m(\ud t,\ud \underline{Z}_t) = \frac{\P([A_t = m]|\UZ_t,\ud \underline{Z}_t)}{\P([A_t = m]|\UZ_t)}.\label{eqn:E_def}
\end{equation}
We expand $E^m(\ud t, \ud \underline{Z}_t)$ as a multivariate
series about $(0,\underline{0})$:
\begin{align}
E^m(\ud t, \ud \underline{Z}_t) &= E^m(0,\underline{0}) + E^m_{\ud t}(0,\underline{0})\ud t + \sum_{j=1}^M E^m_{\ud Z_t^j} (0,\underline{0})\ud Z_t^j\nonumber\\
                                &+ \frac{1}{2}\sum_{j,k=1}^M E^m_{\ud Z_t^j,\ud Z_t^k}(0,\underline{0})\ud Z_t^j \ud Z_t^k + o(\ud t).\label{eqn:multi_series}
\end{align}

By direct evaluation, we obtain:
\begin{align}
&E^m(0,\underline{0}) = 1,\; E^m_{\ud t}(0,\underline{0}) = 0,\nonumber\\
&E^m_{\ud Z_t^j} (0,\underline{0}) = -\frac{1}{\sigma_W^2}\beta_t^j \hat{h}_t,\; j \neq m\nonumber\\
&E^m_{\ud Z_t^m} (0,\underline{0}) = \frac{1}{\sigma_W^2}(1-\beta_t^m) \hat{h}_t,\nonumber\\
&E^m_{\ud Z_t^j,\ud Z_t^j}(0,\underline{0}) = \frac{1}{\sigma_W^4}\beta_t^j(2\beta_t^j-1)\widehat{h^2_t},\; j \neq m \nonumber\\
&E^m_{\ud Z_t^m,\ud Z_t^m}(0,\underline{0}) = \frac{1}{\sigma_W^4}(1-\beta_t^m)(1-2\beta_t^m)\widehat{h^2_t},\nonumber
\end{align}
where $\hat{h}_t := \E[h(X_t)|\UZ_t]$ and
$\widehat{h^2_t}:=\E[h^2(X_t)|\UZ_t]$.

By using It$\hat{\text{o}}$'s rules,
\begin{equation}
\ud Z_t^j \ud Z_t^k = \begin{dcases*}
        \sigma_W^2 \ud t,&\; \text{if} j = k, \\
        0, &\;\text{otherwise.}
        \end{dcases*}\nonumber
\end{equation} This gives
\begin{align}
E^m(\ud t,\ud \underline{Z}_t) &=  1 + \frac{1}{\sigma_W^2}\hat{h}_t\sum_{j=1}^M \beta_t^j(\ud Z_t^m - \ud Z_t^j )\nonumber\\
                               &\quad + \frac{1}{\sigma_W^2}\widehat{h^2_t}  \sum_{j=1}^M \beta_t^j(\beta_t^j - \beta_t^m)\ud t,\label{eqn:Em_approx}
\end{align}
Substituting~\eqref{eqn:Em_approx} to~\eqref{eqn:cal_dbeta} we
otain the expression for $\ud \beta_t^m$
which equals the measurement update part of the continuous-time
filter~\eqref{eqn:filter_for_beta_nonlinear}.
\begin{remark}
During a discrete-time implementation, one can
use~\eqref{eqn:discret_bayes}-\eqref{eqn:Z_t_m} to obtain
association probability. In~\eqref{eqn:discret_bayes}, $L(\ud
Z_t^m)$ is approximated by using particles:
\begin{equation}
L(\ud Z_t^m) \approx \frac{1}{N}\frac{1}{\sqrt{2\pi\sigma_W^2\ud t}} \sum_{i=1}^N \exp\left[-\frac{(\ud Z_t^m - h(X_t^i)\ud t)^2}{2\sigma_W^2\ud t}\right].\nonumber
\end{equation}
\end{remark}

\bibliographystyle{plain}
\bibliography{strings,ACCPF,ACCJPDAFPF,refmtt}

\begin{thebibliography}{10}

\bibitem{Bar-Shalom_IEEE_CSM}
Y.~Bar-Shalom, F.~Daum, and J.~Huang.
\newblock The probabilistic data association filter.
\newblock {\em IEEE Control Systems Magazine}, 29(6):82--100, Dec 2009.

\bibitem{Bar-Shalom_book_88}
Y.~Bar-Shalom and T.~E. Fortmann.
\newblock {\em Tracking and Data Association}.
\newblock Academic Press, San Diego, CA, 1988.

\bibitem{Blackman_book}
S.~S. Blackman.
\newblock {\em Multiple-Target Tracking with Radar Applications}.
\newblock Artech House, Boston, MA, 1986.

\bibitem{BlomBloem2000}
H.~A.~P. Blom and E.~A. Bloem.
\newblock Probabilistic data association avoiding track coalescence.
\newblock {\em IEEE Trans. Automat. Control}, 45(2):247--259, 2000.

\bibitem{BlomBloem2006}
H.~A.~P. Blom and E.~A. Bloem.
\newblock Joint particle filtering of multiple maneuvering targets from
  unassociated measurements.
\newblock {\em Journal of Advancement Information Fusion}, 1:15--36, 2006.

\bibitem{DouFreGor01}
A.~Doucet, N.~de~Freitas, and N.~Gordon.
\newblock {\em Sequential {Monte}-{Carlo} Methods in Practice}.
\newblock Springer-Verlag, April 2001.

\bibitem{gorsalsmi93}
N.~J. Gordon, D.~J. Salmond, and A.~F.~M. Smith.
\newblock Novel approach to nonlinear/non-{Gaussian} {Bayesian} state
  estimation.
\newblock {\em IEE Proceedings F Radar and Signal Processing}, 140(2):107--113,
  1993.

\bibitem{Hue00trackingmultiple}
C.~Hue, J-P.~Le Cadre, and P.~Prez.
\newblock Tracking multiple objects with particle filtering.
\newblock {\em IEEE Trans. Aerospace and Electronic Systems}, 38(3):791--812,
  July 2002.

\bibitem{jazwinski70}
A.~H. Jazwinski.
\newblock {\em Stochastic processes and filtering theory}.
\newblock Academic Press, New York, 1970.

\bibitem{Bar-Shalom_Proc_IEEE}
T.~Kirubarajan and Y.~Bar-Shalom.
\newblock Probabilistic data association techniques for target tracking in
  clutter.
\newblock {\em Proceedings of The IEEE}, 92(3):536--557, 2004.

\bibitem{Herman_2011}
B.~Kragel, S.~Herman, and N.~Roseveare.
\newblock A comparison of methods for estimating track-to-track assignment
  probabilities.
\newblock {\em IEEE Trans. Aerospace and Electronic Systems}, 2011.
\newblock In Press.

\bibitem{Kyriakides_08}
I.~Kyriakides, D.~Morrell, and A.~Papandreou-Suppappola.
\newblock Sequential {M}onte {C}arlo methods for tracking multiple targets with
  deterministic and stochastic constraints.
\newblock {\em IEEE Trans. Signal Process.}, 56(3):937--948, 2008.

\bibitem{Oh_09}
S.~Oh, S.~Russell, and S.~Sastry.
\newblock Markov chain {M}onte {C}arlo data association for multi-target
  tracking.
\newblock {\em IEEE Trans. Automat. Control}, 54(3):481--497, 2009.

\bibitem{Reid79analgorithm}
D.~B. Reid.
\newblock An algorithm for tracking multiple targets.
\newblock {\em IEEE Transactions on Automatic Control}, 24:843--854, 1979.

\bibitem{Ristic_book_2004}
B.~Ristic, S.~Arulampalam, and N.~Gordon.
\newblock {\em Beyond the Kalman Filter: Particle Filters for Tracking
  Applications}.
\newblock Artech House, Boston, MA, 2004.

\bibitem{YangMehtaMeyn11CDC}
T.~Yang, P.~G. Mehta, and S.~P. Meyn.
\newblock Feedback particle filter with mean-field coupling.
\newblock {\em In Proc. of IEEE Conference on Decision and Control}, pages
  7909--7016, December 2011.

\bibitem{YangMehtaMeyn11ACC}
T.~Yang, P.~G. Mehta, and S.~P. Meyn.
\newblock A mean-field control-oriented approach to particle filtering.
\newblock {\em In Proc. of American Control Conference}, pages 2037--2043, June
  2011.

\end{thebibliography}
\end{document}